\documentclass[10pt]{article}

\usepackage{a4wide}
\usepackage{amssymb}
\usepackage{amsfonts}
\usepackage{amsmath}
\input xy
\xyoption{arrow} \xyoption{matrix}

\date{}

\newtheorem{proposition}{Proposition}[section]
\newtheorem{theorem}[proposition]{Theorem}
\newtheorem{lemma}[proposition]{Lemma}

\def\der{\partial }

\def\nFM0{{\nu }_{F,M_0}}
\def\nFN0{{\nu }_{F,N_0}}
\def\nGN0{{\nu }_{G,N_0}}

\def\N0{ {\bf N}_0 }

\def\v{\varphi}
\def\ra{\rightarrow}

\def\Xpm{X^{\pm }}

\def\s{\sigma}

\def\l1{{\lambda}_1}

\def\a{\alpha}
\def\a0{ {\alpha }_0}
\def\a1{ {\alpha }_1}

\def\l{\lambda}
\def\o{\omega}

\def\nFGM0{{\nu }_{F,G,M_0}}

\def\nFN0{{\nu}_{F,N_0}}

\def\sm{{\sigma}^m}

\def\sm1{{\sigma}^{-1}}

\def\smtp1{{\sigma}^{-t+1}}

\def\o{\omega }
\def\S1{S^{-1}}

\def\Xpm1{X^{\pm 1}_1}

\def\sPM1{{\sigma }^{\pm 1}}
\def\sMP1{{\sigma }^{\mp 1 }}

\def\d{\delta}

\def\di{{\rm d.ind}}

\def\L{\Lambda}

\def\CD{{\cal D}}

\def\Ytm1{Y^{t-1}}
\def\Yim1{Y^{i-1}}

\def\CG{{\cal G}}
\def\CH{{\cal H}}

\def\Aut{{\rm Aut}}
\def\bK{\overline{K}}

\def\Der{{\rm Der }}
\def\ad{{\rm ad }}
\def\dim{{\rm dim }}
\def\char{{\rm char }}

\def\ker{ {\rm ker } }

\def\CJ{ {\cal J}}

\def\SL2Z{ {\rm SL}_2({\bf Z}) }

\def\th{ \theta }

\def\Gp1{ G^{1 , 1 } }
\def\P11{ P^{-1 , 1 } }
\def\Pp1{ P^{1 , 1 } }

\def\th{\theta}

\def\nCLsr{{}^\nu\kern-2pt {\cal L}^{\sigma , \rho  }}
\def\nP{{}^\nu \kern-2pt P}
\def\nL{{}^\nu\kern-2pt L}
\def\nLL{{}^\nu\kern-2pt \Lambda}
\def\nPsr{{}^\nu\kern-2pt P^{\sigma , \rho  }}
\def\nLsr{{}^\nu\kern-2pt L^{\sigma , \rho  }}
\def\nuCL{{}^\nu\kern-2pt  {\cal L}}
\def\nCLsr{{}^\nu\kern-2pt {\cal L}^{\sigma , \rho  }}
\def\nCL1m{{}^\nu\kern-2pt {\cal L}^{-1 , 1  }}
\def\x1nu{x^\frac{1}{\nu}}
\def\xm1nu{x^{-\frac{1}{\nu}}}

\def\ra{\rightarrow }

\def\CB{{\cal B}}

\def\CH{ {\cal H}}

\def\nAM0{{\nu }_{{\cal A},M_0}}
\def\nAN0{{\nu }_{{\cal A},N_0}}

\def\End{ {\rm End }}
\def\Der{ {\rm Der }}
\def\CJ{ {\cal J }}

\def\det{ {\rm det }}
\def\ad{ {\rm ad }}

\def\ga{\mathfrak{a}}

\def\SL{{\rm SL}}

\def\di!{\frac{\der^i}{i!}}
\def\dik!{\frac{\der^k_i}{k!}}

\def\id{{\rm id}}

\def\N{\mathbb{N}}

\def\0{\overline{0}}
\def\1{\overline{1}}

\def\Ln1{\L_{n,\overline{1}}}

\def\a1{a_{\overline{1}}}

\def\St{{\rm St}}
\def\S{\Sigma}

\def\grad{{\rm grad}}

\def\vn1{\overrightarrow{n-1}}

\def\Sh{{\rm Sh}}

\def\im{{\rm im}}

\def\Inn{{\rm Inn}}

\def\mJ{\mathbb{J}}
\def\mI{\mathbb{I}}

\def\Cen{{\rm Cen}}

\def\mF{\mathbb{F}}

\def\mT{\mathbb{T}}

\def\mG{\mathbb{G}}
\def\mE{\mathbb{E}}

\def\K1{{\rm K}_1}

\def\hmI1{\widehat{\mI_1}}
\def\tmI1{\widetilde{\mI_1}}
\def\tmJ1{\widetilde{\mJ_1}}
\def\hB1{\widehat{B_1}}
\def\hCB1{\widehat{\CB_1}}

\def\ggu{\mathfrak{u}}

\def\Fix{{\rm Fix}}

\def\UAut{{\rm UAut}}

\def\mJ{\mathbb{J}}

\def\AutLie{ {\rm Aut}_{ {\rm Lie} } }
\def\AutKalg{ {\rm Aut_{K-{\rm alg}}}}
\def\Vir{{\rm Vir}}

\def\mmW1{\mathbb{W}_1}

\def\mQn{\mathbb{Q}_n}
\def\mQ{\mathbb{Q}}
\def\mEn{\mathbb{E}_n}
\def\mE{\mathbb{E}}
\def\Nil{{\rm Nil}}
\def\CX{{\cal X}}
\def\Sym{{\rm Sym}}

\def\divn0{\mathfrak{div}_n^0}
\def\div0mu{\mathfrak{div}_{n, [\mu ]}^0}
\def\din0{\mathfrak{di}_n^0}
\def\ivn0{\mathfrak{iv}_n^0}

\begin{document}

\author{V. V. \  Bavula   
}

\title{The group of automorphisms of the Lie algebra of derivations of a field of rational  functions}

\maketitle

\begin{abstract}

We prove that the group of automorphisms of the Lie algebra $\Der_K (Q_n)$ of derivations of  the field of rational functions  $Q_n=K(x_1, \ldots , x_n)$ over a field of characteristic zero is canonically isomorphic to  the group of automorphisms of the $K$-algebra $Q_n$.

$\noindent $



\end{abstract}


\section{Introduction}

In this paper, module means a left module, $K$ is a
field of characteristic zero and  $K^*$ is its group of units, and the following notation is fixed:
\begin{itemize}
\item $P_n:= K[x_1, \ldots , x_n]=\bigoplus_{\alpha \in \N^n}
Kx^{\alpha}$ is a polynomial algebra over $K$ where
$x^{\alpha}:=x_1^{\alpha_1}\cdots x_n^{\alpha_n}$ and $Q_n:=K(x_1, \ldots , x_n)$ is the field of rational functions,
 \item $G_n:=\AutKalg(P_n)$ and $\mQn :=\AutKalg (Q_n)$;
 \item $\der_1:=\frac{\der}{\der x_1}, \ldots , \der_n:=\frac{\der}{\der
x_n}$ are the partial derivatives ($K$-linear derivations) of
$P_n$,
\item    $D_n:=\Der_K(P_n) =\bigoplus_{i=1}^nP_n\der_i\subseteq E_n:=\Der_K(Q_n) =\bigoplus_{i=1}^nQ_n\der_i$  are  the Lie
algebras of $K$-derivations of $P_n$  and $Q_n$ respectively  where $[\der , \d ]:= \der \d -\d \der $,
\item $\mG_n:=\AutLie(D_n)$ and  $\mEn:=\AutLie (E_n)$,
\item  $\d_1:=\ad (\der_1), \ldots , \d_n:=\ad (\der_n)$ are the inner derivations of the Lie algebras $D_n$ and $E_n$ where $\ad (a)(b):=[a,b]$,
 \item $\CD_n:=\bigoplus_{i=1}^n K\der_i$,
 \item $\CH_n :=\bigoplus_{i=1}^n KH_i$ where $H_1:=x_1\der_1, \ldots , H_n:=x_n\der_n$,
\item  for each natural number $n\geq 2$, $\ggu_n :=
K\der_1+P_1\der_2+\cdots +P_{n-1}\der_n$ is the  {\em Lie algebra of triangular polynomial derivations} (it is a Lie subalgebra of $D_n$) and $ \AutLie(\ggu_n)$ is its group of automorphisms.
    \end{itemize}

\begin{theorem}\label{11Mar13}
\cite{Bav-Aut-Der-Pol}
$\mG_n = G_n$.
\end{theorem}

The aim of the paper is to prove the following theorem.

\begin{theorem}\label{Q11Mar13}
$\mEn = \mQn$.
\end{theorem}

{\it Structure of the proof}. (i) $\mQn\subseteq \mEn$ via the group monomorphism (Lemma \ref{a21Mar13} and (\ref{MMar1}))
$$\mQn\ra \mEn, \;\;  \s \mapsto \s : \der \mapsto \s (\der ):= \s \der \s^{-1}.$$


(ii) Let $\s \in \mEn$. Then $\der_1':=\s (\der_1), \ldots , \der_n':=\s (\der_n)$ are commuting derivations of  $Q_n$  such that  $E_n = \bigoplus_{i=1}^n Q_n \der_i'$ (Lemma \ref{Qc13Mar13}.(2)) and

$\noindent $

(iii) $\bigcap_{i=1}^n \ker_{Q_n}(\der_i') = K$  (Lemma \ref{Qc13Mar13}.(1)).

$\noindent $

(iv)(crux) There exist elements $x_1' , \ldots , x_n'\in Q_n$ such that $\der_i'(x_j')=\d_{ij}$  for  $i,j=1, \ldots , n$ (Lemma \ref{Qc13Mar13}.(3)).

$\noindent $

(v) $\s (x^\alpha \der_i) = x'^\alpha \der_i'$ for all $\alpha \in \N^n$ and $i=1, \ldots , n$ (Lemma \ref{Qc13Mar13}.(6)).

$\noindent $

(vi)  The $K$-algebra homomorphism  $\s' : Q_n \ra Q_n$, $x_i\mapsto x_i'$, $i=1, \ldots , n$ is an automorphism  such that $\s' (q\der_i) = \s' (q) \der_i'$ for all $q\in Q_n$ and $i=1, \ldots , n$.

$\noindent $

(vii) $\Fix_{\mE_n} ( \der_1,\ldots , \der_n, H_1, \ldots , H_n)=\{e\}$ (Proposition \ref{QB11Mar13}.(1)).  Hence,
 $\s = \s'\in \mQ_n$, by (v) and (vi),  i.e. $\mEn = \mQn$.   $\Box $

$\noindent $

{\bf The groups of automorphisms of the Lie algebras $D_n$ and $\ggu_n$}.

\begin{theorem}\label{11Mar13}
 {\rm (Theorem 5.3,  \cite{Bav-Lie-Un-AUT})}  $\AutLie (\ggu_n )\simeq \mT^n\ltimes (\UAut_K(P_n)_n\rtimes( \mF_n' \times  \mE_n ))  $
  where $\mT^n$ is an algebraic  $n$-dimensional torus,   $\UAut_K(P_n)_n$ is an explicit factor group of the group $\UAut_K(P_n)$ of unitriangular polynomial automorphisms, $\mF_n'$ and $\mE_n$ are explicit groups that are isomorphic respectively to the groups $\mI$ and $\mJ^{n-2}$ where
    $\mI := (1+t^2K[[t]], \cdot )\simeq K^{\N}$ and
  $\mJ := (tK[[t]], +)\simeq K^\N$.
\end{theorem}

Comparing the groups $\mG_n$, $\mE_n$ and $\AutLie (\ggu_n)$ we see that the group $\UAut_K(P_n)_n$ of polynomial automorphisms is a {\em tiny} part of the group $\AutLie (\ggu_n)$ but in contrast $\mG_n = G_n$ and  $\mE_n= \mQ_n $.


\begin{theorem}\label{10Mar12}
\cite{Bav-Lie-Un-MON} Every monomorphism of the Lie algebra $\ggu_n$ is an automorphism.
\end{theorem}

Not every epimorphism of the Lie algebra $\ggu_n$ is an automorphism. Moreover, there are countably many distinct ideals $\{ I_{i\o^{n-1}}\, | \, i\geq 0\}$ such that $$I_0=\{0\}\subset I_{\o^{n-1}}\subset I_{2\o^{n-1}}\subset \cdots \subset I_{i\o^{n-1}}\subset \cdots$$ and the Lie algebras $ \ggu_n/I_{i\o^{n-1}}$ and $\ggu_n$ are isomorphic (Theorem 5.1.(1), \cite{Bav-Lie-Un-GEN}).

$\noindent $

{\bf Conjecture}, \cite{Bav-Aut-Der-Pol}. {\em  Every homomorphism of the Lie algebra $D_n$ is an automorphism.}

$\noindent $

The groups of automorphisms of the {\em Witt}  $W_n$ ($n\geq 2$) and the {\em Virasoro}  $\Vir $ Lie algebras were found in \cite{Bav-Aut-Witt-Vir}.



\section{Proof of Theorem \ref{Q11Mar13} }\label{QPPPAAA}

This section can be seen as a proof of Theorem \ref{Q11Mar13}. The proof is split into several statements that reflect `Structure of the proof of Theorem \ref{Q11Mar13}' given in the Introduction.

Let $\CG$ be a Lie algebra and $\CH$ be its Lie subalgebra. The {\em centralizer} $C_\CG (\CH ) := \{ x\in \CG \, | \, [ x, \CH ] =0\}$ of $\CH$ in $\CG$ is a Lie subalgebra of $\CG$. In particular, $Z(\CG ) := C_{\CG }(\CG ) $ is the {\em centre} of the Lie algebra $\CG$. The {\em normalizer} $N_\CG (\CH ) :=\{ x\in \CG \, | \, [ x, \CH ] \subseteq \CH\}$ of $\CH$ in $\CG$ is a Lie subalgebra of $\CG$, it is the largest Lie subalgebra of $\CG$ that contains $\CH $ as an ideal.

Let $V$ be a vector space over $K$. A $K$-linear map $\d : V\ra V$
is called a {\em locally nilpotent map} if $V=\bigcup_{i\geq 1} \ker
(\d^i)$ or, equivalently, for every $v\in V$, $\d^i (v) =0$ for
all $i\gg 1$. When  $\d$ is a locally nilpotent map in $V$ we
also say that $\d$ {\em acts locally nilpotently} on $V$. Every {\em nilpotent} linear map  $\d$, that is $\d^n=0$ for some $n\geq 1$, is a locally nilpotent map but not vice versa, in general.   Let
$\CG$ be a Lie algebra. Each element $a\in \CG$ determines the
derivation  of the Lie algebra $\CG$ by the rule $\ad (a) : \CG
\ra \CG$, $b\mapsto [a,b]$, which is called the {\em inner
derivation} associated with $a$. The set $\Inn (\CG )$ of all the
inner derivations of the Lie algebra $\CG$ is a Lie subalgebra of
the Lie algebra $(\End_K(\CG ), [\cdot , \cdot ])$ where $[f,g]:=
fg-gf$. There is the short exact sequence of Lie algebras
$$ 0\ra Z(\CG ) \ra \CG\stackrel{\ad}{\ra} \Inn (\CG )\ra 0,$$
that is $\Inn (\CG ) \simeq \CG / Z(\CG )$ where $Z(\CG )$ is the {\em centre} of the Lie algebra $\CG$ and $\ad ([a,b]) = [
\ad (a) , \ad (b)]$ for all elements $a, b \in \CG$. An element $a\in \CG$ is called a {\em locally nilpotent element} (respectively, a {\em nilpotent element}) if so is the inner derivation $\ad (a)$ of the Lie algebra $\CG$.

$\noindent $

{\bf The Lie algebra $E_n$}. Since
\begin{equation}\label{EnQnH}
E_n=\bigoplus_{i=1}^n Q_n\der_i= \bigoplus_{i=1}^nQ_n H_i
\end{equation}
every element $\der \in E_n$ is a unique sum $\der =\sum_{i=1}^n a_i\der_i= \sum_{i=1}^n b_iH_i$ where $a_i=x_ib_i\in Q_n$.  The field $Q_n$ is the union $\bigcup_{0\neq f\in P_n} P_{n,f}$ where $P_{n,f}$ is the localization of $P_n$ at the powers of $f$. The algebra $Q_n$ is a localization of $P_{n,f}$. Hence $D_{n,f}:= \Der_K(P_{n,f})= \bigoplus_{i=1}^n P_{n,f}\der_i \subseteq E_n$ and
$$ E_n=\bigcup_{0\neq f\in P_n} D_{n,f}.$$
{\bf $Q_n$ is an $E_n$-module}. The field $Q_n$ is a (left) $E_n$-module: $E_n \times Q_n\ra Q_n$, $(\der, q)\mapsto \der *q$. In more detail, if $\der = \sum_{i=1}^n a_i\der_i$ where $a_i\in Q_n$ then
$$\der * q = \sum_{i=1}^n a_i\frac{\der q}{\der x_i}.$$
The $E_n$-module $Q_n$ is not a simple module since $K$ is an $E_n$-submodule of $Q_n$, and
\begin{equation}\label{QIkerdiK}
\bigcap_{i=1}^n \ker_{Q_n}(\der_i)= K.
\end{equation}

\begin{lemma}\label{a28Mar13}
The $E_n$-module $Q_n/K$  is simple with $\End_{E_n}(Q_n/K)=K\id  $ where $\id$ is the identity map.
\end{lemma}

{\it Proof}. We have to show that for each non-scalar  rational function, say $pq^{-1}\in Q_n$, the $E_n$-submodule $M$ of $Q_n/K$ it  generates coincides  with the $E_n$-module $Q_n/K$. By (\ref{QIkerdiK}), $a_i=\der_i*(pq^{-1})\neq 0$ for some $i$. Then for all elements $u\in Q_n$, $ua_i^{-1}\der_i*(pq^{-1}+K)= u+K$. So, $Q_n/K$ is  a simple $E_n$-module. Let $f\in \End_{E_n}(Q_n/K)$. Then applying $f$ to the equalities $\der_i*(x_1+K)=\d_{i1}$ for $i=1, \ldots , n$, we obtain the equalities
$$ \der_i*f(x_1+K)=\d_{i1} \;\; {\rm for }\;\; i=1, \ldots , n.$$ Hence, $f(x_1+K)\in \bigcap_{i=2}^n \ker_{Q_n/K}(\der_i) \cap \ker_{Q_n/K}(\der_i^2) = (K(x_1)/K)\cap \ker_{Q_n/K}(\der_i^2) =K(x_1+K)$. So, $f(x_1+K) = \l ( x_1+K)$ and so $f=\l \id$, by the simplicity of the $E_n$-module  $Q_n/K$. $\Box $

$\noindent $

{\bf The Cartan subalgebra $\CH_n$ of $E_n$}. A nilpotent Lie subalgebra $C$ of a Lie algebra $\CG$ is called a {\em Cartan subalgebra} of $\CG$ if it coincides with its normalizer. We use often the following obvious observation: {\em An abelian Lie subalgebra that coincides with its centralizer is a maximal abelian Lie subalgebra}.

\begin{lemma}\label{Qa11Mar13}
\begin{enumerate}
\item $\CH_n$ is a Cartan subalgebra of $E_n$.
\item  $\CH_n=C_{E_n}(\CH_n)$ is a maximal abelian Lie subalgebra of $E_n$.
\end{enumerate}
\end{lemma}

{\it Proof}. 2. Clearly, $\CH_n\subseteq C_{E_n}(\CH_n)$. Let $ \der = \sum_{i=1}^n a_iH_i\in C_{E_n}(\CH_n)$ where $a_i\in Q_n$. Then all $a_i\in \cap_{i=1}^n \ker_{Q_n}(H_i) =\cap_{i=1}^n \ker_{Q_n}(\der_i)=K$, by (\ref{QIkerdiK}), and so $\der\in \CH_n$. Therefore, $\CH_n=C_{E_n}(\CH_n)$ is a maximal abelian Lie subalgebra of $E_n$.

1. By statement 2, we have to show that $\CH_n=N:=N_{E_n}(\CH_n)$. Let $\der = \sum_{i=1}^n a_iH_i\in N$, we have to show that all $a_i\in K$. By statement 2, for all $j=1, \ldots , n$, $\CH_n\ni [H_j, \der ] = \sum_{i=1}^nH_j(a_i) H_i$, and so $H_j(a_i)\in K$ for all $i$ and $j$. This condition holds if all $a_i\in K$, i.e. $\der \in \CH_n$. Suppose that $a_i \not\in K$ for some $i$, we seek a contradiction. Then necessarily, $a_i\not\in K(x_1, \ldots , \widehat{x}_j, \ldots , x_n)$ for some $j$. Since $Q_n = K(x_1, \ldots , \widehat{x}_j, \ldots , x_n)(x_j)$, the result follows from the following claim.

{\it Claim: If $ a\in K(x)\backslash K$ then} $H(a) \not\in K$. The field $K(x)$ is a subfield of the series field $K((x)):=\{ \sum_{i>-\infty } \l_i x^i\, | \, \l_i\in K\}$. Since $H(\sum_{i>-\infty } \l_i x^i) = \sum_{i>-\infty } i\l_i x^i$, the Claim is obvious. Then, by the Claim, $H_j(a_i) \not\in K$, a contradiction.  $\Box $

$\noindent $

\begin{lemma}\label{a21Mar13}
\cite{Bav-Aut-Witt-Vir} Let $R$ be a commutative ring such that there exists a derivation $\der \in \Der (R)$ such that $r\der \neq 0$ for all nonzero elements $r\in R$ (eg, $R= P_n,  Q_n$ and $\d = \der_1$). Then the group homomorphism
$$ \Aut (R) \ra \AutLie (\Der (R)), \;\; \s \mapsto \s : \d\mapsto \s (\d ) := \s \d \s^{-1},$$ is a monomorphism.
\end{lemma}


{\bf  The $\mQn$-module $E_n$}. The Lie algebra $E_n$ is a $\mQn$-module,
$$ \mQn\times E_n\ra E_n, \;\; (\s , \der ) \mapsto \s (\der ) := \s \der \s^{-1}.$$
By Lemma \ref{a21Mar13}, the $\mQn$-module $E_n$ is faithful and the map 
\begin{equation}\label{MMar1}
\mQn\ra \mEn, \;\; \s\mapsto \s : \der\mapsto \s (\der ) = \s \der \s^{-1},
\end{equation}
is a group monomorphism. We identify the group $\mQn$ with its image in $\mEn$, $\mQn\subseteq \mEn$. Every automorphism $\s \in \mQn$ is uniquely determined by the elements
$$x_1':=\s (x_1), \; \ldots , \; x_n':=\s (x_n).$$
Let $M_n(Q_n)$ be the algebra of $n\times n$ matrices over  $Q_n$. The matrix  $J(\s) := (J(\s )_{ij}) \in M_n(Q_n)$, where $J(\s )_{ij} =\frac{\der x_j'}{\der x_i}$,   is called the {\em Jacobian matrix} of the automorphism (endomorphism)  $\s$ and its determinant $\CJ (\s ) :=\det \, J(\s)$ is called the {\em Jacobian} of $\s$. So, the $j$'th column of $J(\s )$ is the {\em gradient} $\grad \, x_j':=(\frac{\der x_j'}{\der x_1}, \ldots , \frac{\der x_j'}{\der x_n})^T$  of the polynomial $x_j'$. Then the derivations
$$\der_1':= \s \der_1\s^{-1}, \; \ldots , \; \der_n':= \s\der_n\s^{-1}$$ are the partial derivatives of $Q_n$ with respect to the variables $x_1', \ldots , x_n'$,
\begin{equation}\label{ddp=dxi}
\der_1'=\frac{\der}{\der x_1'}, \; \ldots , \; \der_n'=\frac{\der}{\der x_n'}.
\end{equation}
Every derivation $\der \in E_n$ is a unique sum $\der = \sum_{i=1}^n a_i\der_i$ where $a_i = \der *x_i\in Q_n$. Let  $\der := (\der_1, \ldots , \der_n)^T$ and $ \der' := (\der_1', \ldots , \der_n')^T$ where $T$ stands for the transposition. Then
\begin{equation}\label{dp=Jnd}
\der'=J(\s )^{-1}\der , \;\; {\rm i.e.}\;\; \der_i'=\sum_{j=1}^n (J(\s )^{-1})_{ij} \der_j\;\; {\rm for }\;\; i=1, \ldots , n.
\end{equation}
In more detail, if $\der'=A\der $ where $A= (a_{ij})\in M_n(Q_n)$, i.e. $\der_i=\sum_{j=1}^n a_{ij}\der_j$. Then for all $i,j=1, \ldots , n$,
$$\d_{ij}= \der_i'*x_j'=\sum_{k=1}^na_{ik}\frac{\der x_j'}{\der x_k}$$
where $\d_{ij}$ is the Kronecker delta function. The equalities above can be written in the matrix form as  $AJ(\s) = 1$ where $1$ is the identity matrix. Therefore, $A= J(\s )^{-1}$.

$\noindent $

{\bf The maximal abelian Lie subalgebra $\CD_n$ of $E_n$}.
 Suppose that a group $G$ acts on a set $S$. For a nonempty subset $T$ of $S$, $\St_G(T):=\{ g\in G\, | \, gT=T\}$ is the {\em stabilizer} of the set $T$ in $G$ and $\Fix_G(T):=\{ g\in G\, | \, gt=t$ for all $t\in T\}$  is the {\em fixator} of the set $T$ in $G$. Clearly, $\Fix_G(T)$ is a {\em normal} subgroup of $\St_G(T)$.

\begin{lemma}\label{Qb11Mar13}

\begin{enumerate}
\item $C_{E_n}(\CD_n) =\CD_n$ and so $\CD_n$ is a maximal abelian Lie subalgebra of $E_n$.
\item $\Fix_{\mQn}(\CD_n) = \Fix_{\mQn}(\der_1, \ldots , \der_n) = \Sh_n$.
    \item $\Fix_{\mQn}=(\der_1, \ldots , \der_n , H_1, \ldots , H_n)=\{ e\}$.
        \item $\Cen_{E_n}(\CD_n + \CH_n) = 0$.
\end{enumerate}
\end{lemma}

{\it Proof}. 1. Statement 1 follows from (\ref{QIkerdiK}): Clearly, $\CD_n\subseteq C_{E_n}(\CD_n)$. Let $\der = \sum a_i\der_i\in C_{E_n}(\CD_n)$ where $a_i\in Q_n$. Then all elements $a_i\in \bigcap_{i=1}^n \ker_{Q_n}\der_i= K$, by (\ref{QIkerdiK}), and so $\der \in \CD_n$. So,  $C_{E_n}(\CD_n) =\CD_n$ and as a result  $\CD_n$ is a maximal abelian Lie subalgebra of $E_n$.

2. Let $\s \in \Fix_{\mQn}(\CD_n)$ and  $J(\s ) = (J_{ij})$. By (\ref{dp=Jnd}), $\der = J(\s ) \der $, and so, for all $i,j=1, \ldots , n$,
$\d_{ij} = \der_i*x_j=J_{ij}$,
i.e. $J(\s ) = 1$, or equivalently, by (\ref{QIkerdiK}),
$$x_1'=x_1+\l_1, \ldots , x_n'=x_n+\l_n$$
for some scalars $\l_i\in K$, and so $\s\in \Sh_n$ (since $x_i'-x_i\in \bigcap_{j=1}^n \ker_{Q_n}(\der_j)=K$ for $i=1, \ldots , n$).

3. Let $\s\in \Fix_{\mQn}=(\der_1, \ldots , \der_n , H_1, \ldots , H_n)$. Then $\s \in  \Fix_{\mQn}(\der_1, \ldots , \der_n)=\Sh_n$, by statement 2. So, $\s (x_1) = x_1+\l_1, \ldots , \s (x_n) = x_n+\l_n$ where $\l_i\in K$. Then $x_i\der_i = \s (x_i\der_i) = (x_i+\l_i) \der_i$ for $i=1, \ldots , n$, and so $\l_1=\cdots = \l_n=0$. This means that $\s = e$. So,
 $\Fix_{\mQn}=(\der_1, \ldots , \der_n , H_1, \ldots , H_n)=\{ e\}$.

 4. Statement 4 follows from statement 1 and Lemma \ref{Qa11Mar13}. $\Box $


\begin{lemma}\label{23Mar13}
Let $A$ be a $K$-algebra, $\Der_K(A)$ be the Lie algebra of $K$-derivations of $A$ and $\CD (A)$ be the ring of differential operators on $A$. If the algebra $\CD (A)$ is simple and generated by $A$ and $\Der_K(A)$ then the $\CD (A)$-module $A$ is simple.
\end{lemma}

{\it Proof}. Let $\ga$ be a nonzero $\CD (A)$-submodule  of $A$. So, $\ga$ is an ideal of $A$ such that $ \der (\ga ) \subseteq \ga$ for all $\der \in \Der_K(A)$.  The algebra $\CD :=\CD (A)$ is generated by $A$ and $D$. So, $\CD \ga \subseteq \ga \CD$ and $ \ga \CD \subseteq \CD \ga$, i.e. $\CD \ga = \ga \CD$ is a nonzero ideal of the simple algebra $\CD$. Hence, $1\in \CD \ga$  and so $ 1=\sum_i a_id_i$ for some elements $ d_i\in \CD$ and $a_i\in \ga \subseteq D$. Then
$$ 1=1*1= \sum_i a_id_i*1\in \ga ,$$ hence $\ga = A$, i.e. $A$ is a simple $\CD (A)$-module.  $\Box $





\begin{theorem}\label{Qc11Mar13}

\begin{enumerate}
\item $E_n$ is a simple Lie algebra.
\item $Z(E_n)=\{ 0\}$.
\item $[E_n, E_n]=E_n$.
\end{enumerate}
\end{theorem}

{\it Proof}. 1. (i)  $n=1$, i.e. $E_1=K(x)\der$ {\em is a simple Lie algebra}: We split the proof into several steps.

(a) $D_1:= K[x]\der$ and $W_1:= K[x,x^{-1}]\der$ are simple Lie subalgebras  of $E_1$ (easy).

(b) For all $\l \in K$, $ W_1(\l) :=K[x, (x-\l )^{-1}]$ is a simple Lie subalgebra of $E_1$, by applying the $K$-automorphism $s_\l : x\mapsto x-\l$ of the $K$-algebra $Q_1$ to $W_1$, i.e. $s_\l (W_1) = W_1(\l )$.

(c) {\em For any nonempty subset $I\subset K$, $W_1(I):=W_1(I)_K:=K[x,(x- \l )^{-1} \, | \, \l \in I]\der$ is a simple Lie subalgebra of} $E_1$: Let $\ga$ be a nonzero ideal of $W_1(I)$ and $ 0\neq a\der \in \ga$.  Then either $a\der\in D_1$ or $ 0\neq [ p\der , a\der ] \in D_1\cap \ga $ for some $p\in P_1$. Since $D_1\subseteq W_1(\l )$ for all $\l \in I$ and $W_1(\l )$ are simple Lie algebra, $\ga \cap W_1(\l ) = W_1(\l )$. Hence $\ga = W_1(I)$ since
$$W_1(I) = \bigcup_{\l \in I} W_1(\l ),$$ i.e. $W_1(I)$ is a simple Lie algebra.

(d) If $K$ is an algebraically closed field then $E_1$ is a simple  Lie algebra since $E_1= W_1(K )$.


The algebra $E_1$ is the union $\bigcup_{0\neq f\in P_1} W_1[f^{-1}]$ of the Lie algebras $ W_1[f^{-1}]:= P_{1,f}\der$ where $P_{1,f}$ is the localization of $P_1$ at the powers of the element $f$. Let $ \ga$ be the ideal of $E_1$ generated by  a nonzero element $ a=pq^{-1}\der$ for some $pq^{-1}\in Q_1$. Clearly, $a\in W_1[(fq)^{-1}]$ for all nonzero elements $f\in P_1$ and  $E_1= \bigcup_{0\neq f\in P_1}W_1[(fg)^{-1}]$. So, to finish the proof of (i) it suffices to show that all the algebras $W_1[f^{-1}]$ are simple.

(e) $A:= W_1[f^{-1}]$ {\em is a simple Lie algebra for all} $0\neq f\in P_1$: Let $K':=K(\nu_1, \ldots , \nu_s)$ be the subfield of the algebraic closure $\bK$ of $K$ generated by the roots $\nu_1, \ldots , \nu_s$ of the polynomial $f$ and $G= {\rm Gal} (K'/K)$ be the Galois group of the finite Galois field extension $K'/K$ (since $\char (K)=0$). Let $K'=\oplus_{i=1}^d K\th_i$ for some elements $\th_i\in K'$ and $\th_1=1$. By (c),
$$A':= K'[x,f^{-1}]\der = W_1(\nu_1, \ldots , \nu_s)_{K'}$$ is a simple Lie $K'$-algebra. Let $ a\in A\backslash \{ 0\}$, $\ga $ and d $\ga'$ be the ideals in $A$ and $A'$ respectively that are generated by the element $a$. Then $\ga' = A'$, by (c). Notice that $A' = \sum_{i=1}^d\th_i A$ and for $a'= \sum_{i=1}^d \th_ia_i, b=\sum_{i=1}^d \th_ib_i\in A'$ where $a_i, b_i\in A$, $[a',b]=\sum_{i=1}^d \th_i\th_j[ a_i, b_j]$. Moreover,  every element in $A' = \ga'$ is a linear combination of several commutators in $A'$ (where $c= \sum_{i=1}^d \th_kc_k\in A'$ and $c_k\in A$),
\begin{equation}\label{abcs}
[a,[a',  \ldots [b,c]\ldots ] = \sum \th_i\cdots \th_j\th_k [a,[a_i, \ldots [b_j,c_k]\ldots ].
\end{equation}
The {\em symmetrization map} $\Sym : K'\ra K$, $\l \mapsto |G|^{-1} \sum_{g\in G} g(\l )$, is a surjection such that $\Sym (\mu ) = \mu$ for all $\mu \in K$. Clearly, $K'(x)/ K(x)$ is a Galois field extension with the Galois group $G$ where the elements of $G$ act trivially on the element $x$. So, the symmetrization map $\Sym$ can be extended to the surjection $ K'(x) \ra K(x)$ by the same rule, and then to the surjection
 $A'\ra A$, $f\der \mapsto \Sym (f) \der $.

 Each element $e\in A\subseteq A'$, can be expressed as a finite sum of elements in  (\ref{abcs}). Then applying $\Sym$, we see that $e$ is a linear combination of elements (commutators) from $\ga$, i.e. $A$ is a simple Lie algebra.

 (ii) $E_n$ {\em is a simple Lie algebra for } $n\geq 2$: Let $a\in E_n\backslash \{ 0\}$ and $\ga = (a)$ be the ideal in $E_n$ generated by the element $a= \sum_{i=1}^na_i\der_i$ where $a_i\in Q_n$.

 (a) $\ga \cap D_n\neq 0$: If $a\in D_n$ then there is nothing to prove. Suppose that $a\not\in D_n$.

 (a1) Suppose that $a_i\in K(x_i)$ for all $i$. Then $a_i\not\in K[x_i]$ for some $i$ (since $a\not\in D_n$), and so
 $$ \ga \ni [ H_i , a] = H_i(a_i) \der_i \in K(x_i) \der_i \backslash \{ 0\}.$$
By (i), $\der_1\in \ga \cap D_n$.

(a2) Suppose that $a_i\not\in K(x_i)$ for some $i$. Then $\der_j(a_i) \neq  0$ for some $j\neq i$. Let $q\in P_n$ be the common denominator of the fractions $a_1, \ldots , a_n$, that is $a_1=p_1q^{-1}, \ldots , a_n = p_nq^{-1}$ for some elements $p_i\in P_n$. For all $n\geq 2$,
$$ D_n\cap \ga \ni [ q^n\der_j , a] = q^n\der_j (a_i) \der_i + \sum_{k\neq i } (\ldots ) \der_k\neq 0.$$
(b) $\ga = D_n$ since $D_n$ is a simple Lie algebra, \cite{Bav-Aut-Der-Pol}.

(c) $\ga \supseteq K(x_i)\der_i$ {\em for} $i=1, \ldots , n$: In view of symmetry it suffices to prove that $\ga \supseteq K(x_1)\der_1$. Notice that for all $u\in Q_n$ and $i=2, \ldots , n$,
$$ \ga \ni [ u\der_1, x_1\der_i] = u\der_i-x_1\der_i(u) \der_1.$$ Therefore, $\ga + Q_n\der_1= E_n$. The field of rational functions $Q_n = Q_n(K)$ can be seen as the field of rational functions $Q_n(K) = Q_{n-1}(K')$ where $K' = K(x_1)$. Then $$E_{n-1}':= \Der_{K'}(Q_{n-1}(K'))= \bigoplus_{i=2}^n Q_{n-1}(K') \der_i= \bigoplus_{i=2}^n Q_n \der_i.$$
By Lemma \ref{23Mar13}, the $E_{n-1}'$-module $Q_{n-1}'/K'= Q_n/K(x_1)$ is simple. The Lie algebra $E_{n-1}'$ is a Lie subalgebra of $E_n$, and $E_n$ can be seen as a left $E_{n-1}'$-module  with respect to the adjoint action. The ideal $\ga$ of $E_n$ is an $E_{n-1}'$-submodule of $E_n$. The Lie algebra $K(x_1)\der_1$ is simple and $\ga \cap K(x_1)\der_1$ is a nonzero ideal of it (by (b)). Therefore, $K(x_1)\der_1\subseteq \ga $. The $E_{n-1}'$-module $E_n /\ga = (\ga +Q_n\der_1) / \ga \simeq Q_n\der_1/ \ga \cap Q_n\der_1$ is an epimorphic image of the simple $E_{n-1}'$-module $Q_n / K(x_1)$ via
$$\v : Q_n/ K(x_1)\ra Q_n\der_1/ \ga \cap Q_n\der_1, \;\; u+K(x_1)\mapsto u\der_1+\ga \cap Q_n\der_1, $$
with $0\neq (P_n+K(x_1))/ K(x_1)\subseteq \ker (\v )$. Therefore, $Q_n\der_1 = \ga \cap Q_n\der_1\subseteq \ga$, and so $E_n = \ga +Q_n\der_1 = \ga$. So, $E_n$ is a simple Lie algebra.

2 and 3. Statements 2 and 3 follow from statement 1.
 $\Box $


\begin{lemma}\label{b23Mar13}
For all nonzero elements $q\in Q_n$ and $i=1, \ldots , n$, $C_{E_n}(qP_n\der_i) = \{ 0\}$.
\end{lemma}

{\it Proof}. Let $c\in C_{E_n}(qP_n\der_i)$. Then for all elements $p\in P_n$,
$$ 0=[c,qp\der_i]= c(p)\cdot q\der_i +p[c,q\der_i]= c(p)\cdot q\der_i.$$ Then $c(p)=0$ for all $p\in P_n$, and so $c=0$. $\Box $


\begin{proposition}\label{B11Mar13}
\cite{Bav-Aut-Der-Pol}\; $\Fix_{\mG_n} (\der_1, \ldots , \der_n, H_1, \ldots , H_n) = \{ e\}$.
\end{proposition}

Let $d_1, \ldots , d_n$ be a commuting linear maps acting in a vector space $E$. Let $\Nil_E(d_1, \ldots , d_n) := \{ e\in E\, | \, d_i^je=0$ for all $i=1, \ldots , n$ and some $j=j(e)\}$. Let  $\Nil_{E_n}(\CD_n):= \Nil_{E_n}(\d_1, \ldots , \d_n)$. Clearly, $\Nil_{E_n}(\CD_n)=D_n$ is a Lie subalgebra of $E_n$.

\begin{proposition}\label{QB11Mar13}

\begin{enumerate}
\item $\Fix_{\mEn} (\der_1, \ldots , \der_n, H_1, \ldots , H_n) = \{ e\}$.
\item $\Fix_{\mEn} (\der_1, \ldots , \der_n) = \Sh_n$.
\end{enumerate}
\end{proposition}

{\it Proof}. 1. Let $\s\in F:= \Fix_{\mEn} (\der_1, \ldots , \der_n, H_1, \ldots , H_n)$. We have to show that $\s = e$. Then $\s^{-1}\in F$ and $\s^{\pm 1}(\Nil_{E_n} (\CD_n))\subseteq \Nil_{E_n} (\CD_n)$, i.e. $\s (D_n)= D_n$ since $\Nil_{E_n} (\CD_n) = D_n$. So, $\s |_{D_n}\in
 \Fix_{\mG_n} (\der_1, \ldots , \der_n, H_1, \ldots , H_n)=\{ e\}$  (Proposition  \ref{B11Mar13}), i.e. $\s (\der ) = \der $ for all $\der \in D_n$. Let $0\neq \d \in E_n$. Then $\d = q^{-1} \der$ for some $0\neq q\in P_n$ and $\der \in D_n$. Now, $[q^2p\der_i, \d ] = \der'\in D_n$ for all $p\in  P_n$. Applying $\s$ to the equality yields the equality $[q^2p\der_i , \s (\d )] = \der'$. By taking the difference, we obtain
 $ \s (\d ) - \d \in C_{E_n} (q^2P_n \der_i) = \{ 0\}$, by Lemma \ref{b23Mar13}, hence $\s =e$.

2. Clearly, $\Sh_n \subseteq F:=\Fix_{\mEn} (\der_1, \ldots , \der_n)$. Let $\s \in F$ and $H_i':=\s (H_i), \ldots , H_n':=\s (H_n)$. Applying the automorphism $\s$ to the commutation relations $[\der_i, H_j]=\d_{ij}\der_i$ gives the relations $[\der_i, H_j']=\d_{ij}\der_i$. By taking the difference, we see that $[\der_i, H_j'-H_j]=0$ for all $i$ and $j$. Therefore, $H_i'=H_i+d_i$ for some elements  $d_i\in C_{E_n}(\CD_n) = \CD_n$ (Lemma \ref{Qb11Mar13}.(1)),  and so $d_i=\sum_{j=1}^n \l_{ij}\der_j$ for some elements $\l_{ij}\in K$. The elements $H_1', \ldots , H_n'$ commute, hence
$$ [ H_j, d_i]= [H_i, d_j] \;\; {\rm for \; all}\;\; i,j,  $$
or equivalently,
$$ \l_{ij}\der_j= \l_{ji}\der_i \;\; {\rm for \; all}\;\; i,j.  $$
This means that $\l_{ij}=0$ for all $i\neq j$, i.e. $$H_i'= H_i+\l_{ii}\der_i= (x_i+\l_{ii})\der_i= s_\l (H_i)$$
where $s_\l \in \Sh_n$, $s_\l (x_i) = x_i+\l_{ii}$ for all $i$. Then $s_\l^{-1}\s \in \Fix_{\mEn} (\der_1, \ldots , \der_n, H_1, \ldots , H_n) = \{ e\}$  (statement 1), and so $\s = s_\l \in \Sh_n$. $\Box $

$\noindent $

{\bf The automorphism $\nu$}. Let $\nu$ be the $K$-automorphism of $Q_n$ given by the rule $\nu (x_i) = x_i^{-1}$ for $i=1, \ldots , n$. Then
\begin{equation}\label{nua}
 \nu (\der_i) = -x_iH_i, \;\; \nu (H_i) = -H_i, \;\; \nu (x_iH_i) = -\der_i, \;\; i=1, \ldots , n.
\end{equation}
By (\ref{nua}), the elements $X_1:=x_1H_1, \ldots , X_n:= x_nH_n$ commute and the next lemma follows from Lemma \ref{Qb11Mar13} and Proposition \ref{QB11Mar13} since $\CX_n:=\nu (\CD_n) = \bigoplus_{i=1}^n KX_i$.

\begin{lemma}\label{x23Mar13}
\begin{enumerate}
\item $C_{E_n}(\CX_n) = \CX_n $ is a maximal abelian Lie subalgebra of $E_n$.
    \item $\Fix_{\mQn} (X_1, \ldots , X_n) =\Fix_{\mEn} (X_1, \ldots , X_n) = \Sh_n$.
\item $\Fix_{\mQn} (X_1, \ldots , X_n, H_1, \ldots , H_n) = \Fix_{\mEn} (X_1, \ldots , X_n, H_1, \ldots , H_n) =\{ e\}$.
    \end{enumerate}
\end{lemma}

The following lemma is well-known and it is easy to prove.

\begin{lemma}\label{Aslice}
Let $\der$ be a locally nilpotent derivation of a commutative $K$-algebra $A$ such that $\der (x) =1$ for some element $x\in A$. Then $A= A^\der [x]$ is a polynomial algebra over the ring $A^\der := \ker (\der )$ of constants of the derivation $\der$ in the variable $x$.
\end{lemma}





The next lemma is the core of the proof of Theorem \ref{Q11Mar13}.
\begin{lemma}\label{Qc13Mar13}
Let $\s \in \mEn$,  $\der_1':=\s (\der_1), \ldots , \der_n':= \s (\der_n)$ and $\d_1':=\ad (\der_1'), \ldots , \d_n':= \ad (\der_n')$. Then
\begin{enumerate}
\item $\der_1', \ldots , \der_n'$ are commuting derivations of $Q_n$ such that $\bigcap_{i=1}^n\ker_{Q_n}(\der_i')=K$.
\item  $E_n = \bigoplus_{i=1}^n Q_n \der_i'$.
\item For each $i=1, \ldots , n$, $\s (x_i\der_i) = x_i'\der_i'$ for some elements $x_i'\in Q_n$. The elements $x_1', \ldots , x_n'$ are algebraically independent and  $\der_i'(x_j') = \d_{ij}$ for $i,j=1, \ldots , n$.
\item $\Nil_{Q_n} (\der_1', \ldots , \der_n') = P_n'$ where $P_n' := K[x_1', \ldots , x_n']$.
\item $\Nil_{E_n} (\d_1', \ldots , \d_n') = \bigoplus_{i=1}^n P_n'\der_i'$.
\item $\s (x^\alpha \der_i) = x'^\alpha \der_i'$ for all $\alpha \in \N^n$ and $i=1, \ldots , n$.
\item  $\s': Q_n\ra Q_n$, $x_i\mapsto x_i'$, $ i=1, \ldots , n$ is a $K$-algebra homomorphism (statement 3) such that $ \s' (a\der_i) = \s' (a) \s (\der_i)$.
\item The $K$-algebra homomorphism $\s'$ is an automorphism.
\end{enumerate}
\end{lemma}

{\it Proof}. 1.  The elements $\der_1, \ldots , \der_n$ are commuting derivations, hence so are  $\der_1', \ldots , \der_n'$.   Let $\l \in \bigcap_{i=1}^n \ker_{Q_n}(\der_i')$. Then
 $$\l \der_1' \in C_{E_n}(\der_1', \ldots , \der_n')=\s (C_{E_n}(\der_1, \ldots , \der_n))= \s (C_{E_n}(\CD_n))= \s (\CD_n) = \s (\bigoplus_{i=1}^n K\der_i) = \bigoplus_{i=1}^n K\der_i', $$
since $C_{E_n}(\CD_n)=\CD_n$, Lemma \ref{Qb11Mar13}.(1). Then $\l \in K$ since otherwise the infinite dimensional space $\bigoplus_{i\geq 0} K\l^i \der_1'$ would be a subspace of the  finite dimensional space $\s (\CD_n)$.

2.  It suffices to show that the elements $\der_1', \ldots , \der_n'$ of the $n$-dimensional (left) vector space $E_n$ over the field $Q_n$ are $Q_n$-linearly independent (the key reason for that is statement 1). Let $V=\sum_{i=1}^nQ_n\der_i'$. Suppose that $m:=\dim_{Q_n}(V)<n$, we seek a contradiction. Up to order, let $\der_1', \ldots , \der_m'$ be a $Q_n$-basis of $V$. Then
$\der_{m+1}=\sum_{i=1}^m a_i\der_i'$ for some elements $a_i\in Q_n$. By applying $\d_j'$ ($j=1, \ldots , n$), we see that $0=\sum_{i=1}^m \der_j'(a) \der_i'$, and so $a_i\in \bigcap_{i=1}^n \ker_{Q_n}(\der_j') = K$, by statement 1. This means that the elements $\der_1', \ldots , \der_m'$ are $K$-linearly dependent, a contradiction.

3. Let $H_i':= \s (x_i\der_i)$ for $i=1, \ldots , n$. By statement 2, $H_i'=\sum_{i=1}^n a_{ij}\der_j'$ for some elements $a_{ij}\in Q_n$. Applying the automorphism $\s$ to the relations  $ \d_{ij} \der_j= [ \der_j, H_i]$ yields the relations
$$\d_{ij} \der_i'= \sum_{i=1}^n \der_j'(a_{ik})\der_k'.$$
 Let $x_i':= a_{ii}$. Then $ \der_j'(x_i') = \d_{ji}$ and $\der_j'(a_{ik})=0$ for all $k\neq i$. By statement 1, $a_{ik}\in K$ for all $i\neq k$. Now,
 $$H_i':= x_i'\der_i'+\sum_{j\neq i} a_{ij}\der_j'.$$ The elements $H_1', \ldots , H_n'$ commute, hence for all $i\neq j$, $0=[H_i', H_j']= - a_{ji}\der_i' +a_{ij} \der_j'$, and so $a_{ij}=0$. Therefore, $H_i' = x_i'\der_i'$.

The equalities $\der_i'(x_j')=\d_{ij}$ imply that the elements $x_1',\ldots , x_n'\in Q_n$ are algebraically independent over $K$: Suppose that $f(x_1',\ldots , x_n')=0$ for some nonzero polynomial $f(t_1, \ldots , x_n) \in K[t_1, \ldots , x_n]$.  We can assume that the (total) degree $\deg (f)$ is the least possible. Clearly, $f\not\in K$, hence $\frac{\der f}{\der x_i}\neq 0$ for some $i$ and $\deg (\frac{\der f}{\der x_i})<\deg (f)$, but $\frac{\der f}{\der x_i} (x_1',\ldots , x_n')= \der_i (f(x_1', \ldots , x_n'))=\der_i(0)=0$, a contradiction.

4. Let $\CD_n' = \sum_{i=1}^n K\der_i'$ and $N= \Nil_{Q_n}(\CD_n')$. By statement 3 and Lemma \ref{Aslice}, $$N=N^{\CD_n'} [ x_1' ,\ldots , x_n']= K [ x_1' ,\ldots , x_n']$$ since $K\subseteq N^{\CD_n'}\subseteq Q_n^{\CD_n'}=K$ (by statement 1).

5. Let $\der= \sum_{i=1}^n a_i\der_i'\in N:= \Nil_{E_n} (\d_1', \ldots , \d_n')$ where $a_i\in Q_n$ (statement 2). For all $\alpha \in \N^n$,
$$ \d'^\alpha (\der ) = \sum_{i=1}^n \der'^\alpha (a_i) \der_i'$$ where $ \d'^\alpha := \prod_{i=1}^n \d_i'^{\alpha_i}$, $\d_i' = \ad (\der_i')$ and $ \der'^\alpha := \prod_{i=1}^n \der_i'^{\alpha_i}$. So, $\d'^\alpha (a_i) =0$  iff  $\der'^\alpha(a_i) =0$ for $i=1, \ldots , n$ (statement 2). Now, statement 5 follows from statement 4.

6. First, let us show that, by  induction on $|\alpha |$, that $\s (x^\alpha \der_i)-x'^\alpha\der_i'\in \Cen_{E_n}(\CD_n')=\CD_n'$  (Lemma \ref{Qb11Mar13}.(1)). The initial case when $|\alpha | =0$ is obvious. So, let $|\alpha | >0$. Then
\begin{eqnarray*}
 [\der_j', \s ( x^\alpha \der_i) - x'^\alpha \der_i'] &=& \s ([\der_j, x^\alpha \der_i]) - \alpha_jx'^{\alpha - e_j} \der_i' = \s (\alpha_j x^{\alpha - e_j} \der_i) - \alpha_jx'^{\alpha - e_j} \der_i'\\
  &=& \alpha_jx'^{\alpha - e_j} \der_i'-\alpha_jx'^{\alpha - e_j} \der_i'=0. \end{eqnarray*}
Therefore, $\s (x^\alpha \der_i) = x'^\alpha\der_i' +\sum \l_{ij} \der_j'$ for some scalars $\l_{ij} = \l_{ij}(\alpha ) \in K$. Notice that $$\s (H_i) = \s (x_i\der_i) = x_i' \der_i':= H_i',$$ by the definition of the elements $x_i'$. Since $|\alpha |>0$, $\alpha_j\neq 0$ for some $j$. Applying the automorphism $\s $ to the equalities $(\alpha_j-\d_{ij})  x^\alpha \der_i = [ H_j, x^\alpha \der_i]$ we have (we may assume that $x^\alpha \der_i \neq H_i$)
\begin{eqnarray*}
 (\alpha_j-\d_{ij}) (x'^\alpha\der_i' +\sum_{k=1}^n \l_{ik} \der_k') &=& \s ((\alpha_j-\d_{ij}) x^\alpha\der_i) = \s ([  H_j, x^\alpha \der_i])=[H_j', x'^\alpha \der_i' +\sum_{k=1}^n\l_{ik} \der_k']\\
  &=&  (\alpha_j -\d_{ij}) x'^\alpha\der_i' -\l_{ij} \der_j',
\end{eqnarray*}
and so $(\alpha_j-\d_{ij} +1) \l_{ij}=0$ and $ (\alpha_j -\d_{ij}) \l_{ik}=0$ for all $k\neq j$. This means that all $\l_{is}=0$.

7.  By statement 3, $\s'$ is a $K$-algebra homomorphism such that $\im (\s') =Q_n':=K(x_1', \ldots , x_n')$. By statement 3, for all elements $a\in Q_n$,
$$ \der_i'\s'(a) = \s'\der_i (a)$$ since $\der_i'$ acts as $\frac{\der}{\der x_i'}$ on $Q_n'$.

 Let $a=pq^{-1}\neq 0 $ where $p,q\in P_n$. Then, for all $r\in q^2P_n$, $[ a\der_i, r\der_i]=(a\der_i(r)- \der_i(a) r)\der_i\in P_n\der_i$. By applying $\s$, we have the equality
$$[ \s (a\der_i) , \s' (r) \der_i'] = \s' ( a\der_i(r)-\der_i(a) r)\der_i'.$$
On the other hand,
\begin{eqnarray*}
[\s' (a)\der_i', \s' (r) \der_i']&= & (\s' (a) \der_i'\s' (r) - \der_i' \s' (a) \s' (r)) \der_i'= (\s' (a) \s' \der_i(r) - \s' \der_i(a)\s' (r)) \der_i'   \\
 &=& \s' ( a\der_i(r) - \der_i(a) r)\der_i'.
\end{eqnarray*}
Hence,
$$
 \s (a\der_i) -\s' (a) \der_i' \in  C_{E_n} (\s' (q^2P_n)\der_i') = C_{E_n}(\s (q^2P_n\der_i)) = \s (C_{E_n} (q^2P_n\der_i)) = \s (C_{E_n}(q^2P_n\der_i)) =0,
 $$
 by Lemma \ref{b23Mar13}.
Therefore, $\s (a\der_i) = \s' (a)\s (\der_i)$.

8. Since $ \s (Q_n\der_i) = \s' (Q_n) \der_i'$ for all $i=1, \ldots , n$ (statement 7), we must have $\s' (Q_n) = Q_n$, by statement 2, and so $\s'\in \mQ_n$. $\Box$

$\noindent $

{\bf Proof of Theorem \ref{Q11Mar13}}. Let $\s \in \mE_n$. By Corollary \ref{Qc13Mar13}.(8), we have the automorphism
 $\s'\in \mQ_n$ such that, by Lemma \ref{Qc13Mar13}.(3,6), $\s'^{-1}\s \in \Fix_{\mE_n}(\der_1, \ldots, \der_n , H_1, \ldots , H_n)=\{ e\}$ (Proposition \ref{QB11Mar13}).
 Therefore, $\s = \s'$ and so  $\mE_n = \mQ_n$. $\Box$


$${\bf Acknowledgements}$$

 The work is partly supported by  the Royal Society  and EPSRC.

\small{

Department of Pure Mathematics

University of Sheffield

Hicks Building

Sheffield S3 7RH

UK

email: v.bavula@sheffield.ac.uk}


\begin{thebibliography}{99}









\bibitem{Bav-Lie-Un-MON}  V. V. Bavula, Every monomorphism of the Lie algebra of triangular polynomial derivations is an automorphism, {\em C. R. Acad. Sci. Paris, Ser. I}, {\bf 350} (2012) no. 11-12,  553--556.  (Arxiv:math.AG:1205.0797).

\bibitem{Bav-Lie-Un-GEN} V. V. Bavula,  Lie algebras of  triangular polynomial derivations and an
isomorphism criterion for their Lie factor algebras,  {\it
Izvestiya: Mathematics}, (2013), in print. (Arxiv:math.RA:1204.4908).






\bibitem{Bav-Lie-Un-AUT} V. V. Bavula,  The groups of automorphisms of the Lie algebras of triangular polynomial derivations, Arxiv:math.AG/1204.4910.

\bibitem{Bav-Aut-Der-Pol}  V. V. Bavula, The group of automorphisms of the Lie algebra of derivations of a polynomial algebra. Arxiv:math.RA:1304.6524.

\bibitem{Bav-Aut-Witt-Vir} V. V. Bavula, The groups of automorphisms of the Witt $W_n$ and Virasoro Lie algebras. Arxiv:math.RA:1304.6578.
























\end{thebibliography}
\end{document}